\documentclass[11pt,reqno]{amsart}

\usepackage{latexsym}
\usepackage{amsfonts}
\usepackage{amsmath}
\usepackage{url}
\usepackage{graphicx}
\usepackage{color}

\usepackage{epsfig,psfrag,verbatim,amsfonts}

\usepackage{amssymb,amsmath}

\def\bl{\begin{lemma}}
\def\el{\end{lemma}}
\def\bth{\begin{theorem}}
\def\eth{\end{theorem}}
\def\bc{\begin{corollary}}
\def\ec{\end{corollary}}
\def\bcj{\begin{conjecture}}
\def\ecj{\end{conjecture}}
\def\bpr{\begin{proposition}}
\def\epr{\end{proposition}}
\def\bde{\begin{definition}}
\def\ede{\end{definition}}

\newcommand{\be}{\begin{eqnarray}}
\newcommand{\ee}{\end{eqnarray}}

\renewcommand{\and}{\hbox{ {\rm and} }}

\newtheorem{theorem}{Theorem}[section]
\newtheorem{definition}{Definition}[section]
\newtheorem{lemma}[theorem]{Lemma}

\newtheorem{corollary}[theorem]{Corollary}
\newtheorem{proposition}[theorem]{Proposition}
\newtheorem{conjecture}[theorem]{Conjecture}

\theoremstyle{definition}
\numberwithin{equation}{section}

\begin{document}
\title{Percolation and coarse conformal uniformization}
\author{Itai Benjamini }

\date{June 2015}

\maketitle

\begin{abstract}
We formulate conjectures regarding percolation on planar triangulations suggested  by assuming (quasi) invariance under coarse conformal uniformization.
\end{abstract}

\section{Introduction}

We start with a quick panoramic  overview.
A {\em conformal map}, between planar domains, is a function that infinitesimally preserves angles.
The derivative  of a conformal map  is everywhere a scalar times a rotation.

{\em Riemann's mapping theorem} states that any open simply connected domain of the Euclidean plane admits a bijective conformal map to the open unit disk.
In the 1940's Shizuo Kakutani observed that two dimensional Brownian motion  is conformal invariant, up to a time reparametrization. Therefore the scaling limit
of simple random walks on the Euclidean grid is conformal invariant.

In 2000 Stas Smirnov \cite{Sm}  proved that the scaling limit of critical Bernoulli site percolation on the triangular lattice is conformal invariant.

Poincar\'{e} (1907) proved that every simply connected Riemann surface is conformally equivalent to one of the following three surfaces:
the open unit disk, the Euclidean plane, or the Riemann sphere. In particular it admits a Riemannian metric of constant curvature.
This classifies Riemannian surfaces as elliptic (the shpere), parabolic (Euclidean), and hyperbolic (negatively curved).

\medskip
 The uniformization theorem is a generalization of the Riemann mapping theorem from proper simply connected open subsets of the plane to arbitrary simply connected Riemannian  surfaces.

Conformal invariance of Brownian motion extends to the context of the uniformization.
A  simply connected Riemann surface is conformally equivalent to the hyperbolic plane iff the Brownian motion is transient.
\medskip

{\em How does surface uniformiztion  manifest  itself in the context of percolation?}
\medskip

Below we suggest  that   the discrete  setup of planar triangulations is  natural for this problem.

Recall, every planar graph admits a circle packing,  Koebe (1936).
\medskip

In 1995 He and Schramm \cite{HS} proved a discrete uniformization theorem for triangulations:  Let $G$ be the $1$-skeleton of a
triangulation of an open disk. If the random walk on $G$ is
{\em recurrent}, then $G$ is circled packed in the Euclidean plane.
Conversely, if the degrees of the vertices in $G$ are
bounded and the random walk on $G$ is {\em transient}, then $G$ is circle packed in the unit disc.

\medskip

For   an extended version of the Brooks, smith, Stone and Tutte (1940) square tiling theorem,
a  related discrete uniformization theorem for graphs using squares, see  \cite{BS0}.
\medskip

\begin{figure}[t]
\centering
\includegraphics[trim = 2em .5em 2em .5em, clip, height=.45\textwidth]{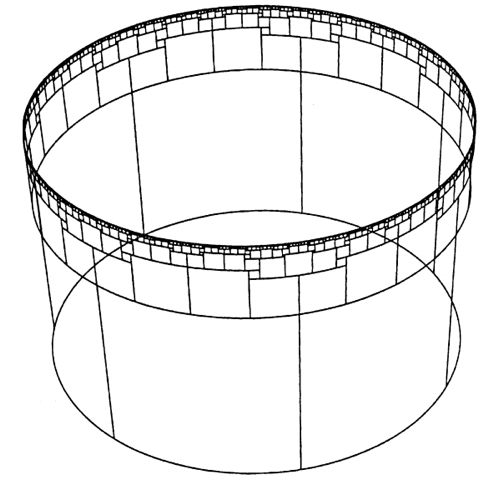} \qquad \includegraphics[trim = 1em 10.5em 3em 10em, clip,height=.45\textwidth]{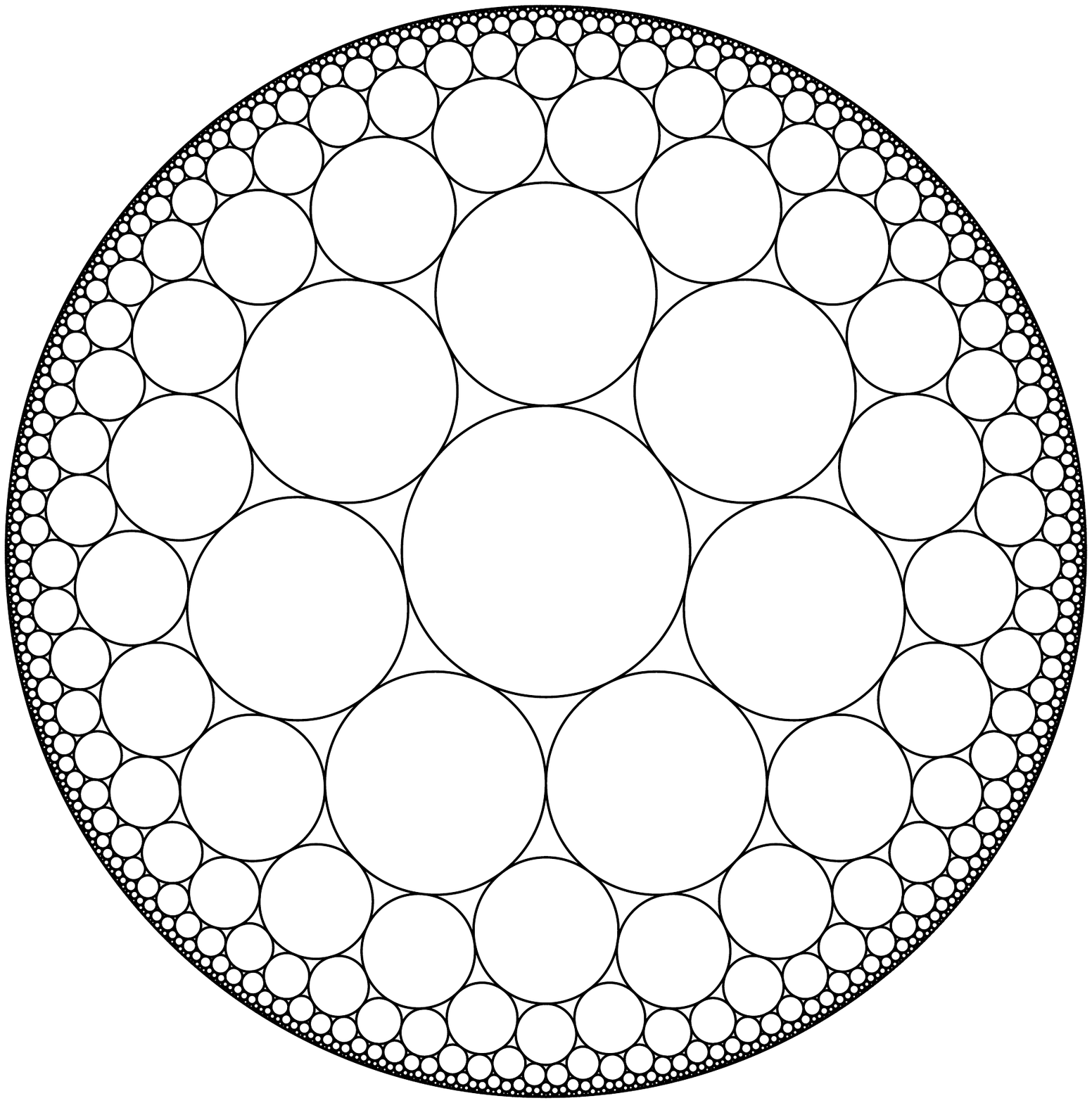}
\caption{The square tiling and the circle packing of the 7-regular  triangulation.}
\end{figure}

In $p$-Bernoulli site percolation, each vertex is declared open independently with probability $p$, and clusters are connected components of open vertices.
A graph is transient if simple random on the graph returns to the origin finitely many times almost surely and otherwise is called recurrent.

\section{Two conjectures}

Let $G$ be the $1$-skeleton of a bounded degree triangulation of an open disk.

\begin{conjecture}\label{c1}
Assume $G$ is transient, then $1/2$-Bernoulli site percolation on $G$ admits an infinite cluster a.s.
\end{conjecture}

We don't know it even  for any  fixed $p >1/2$.

The motivation for the conjecture is outlined below and is based on conformal invariance of percolation.
After more than two decades of thorough research,
conformal invariance of critical Bernoulli percolation was established essentially only for the triangular lattice \cite{Sm}.

One reason to be slightly skeptical about the conjecture is that for critical percolation on the triangular lattice,
the probability the cluster of the origin reaches distance $r$ decays polynomially in $r$ \cite{SW},
while there are transient triangulations of volume  growth $r^2 \log^3 r$.
\medskip

\noindent
{\em A heuristic}

Tile the unit square with (possibly infinity number) of squares of varying sizes so
that at most three squares meet at corners. Color each square black or white with
equal probability independently.
\begin{conjecture}\label{c2}
Show that there is a universal $c > 0$, so that the probability of a
black left right crossing is bigger than $c$.
\end{conjecture}

If true, the same should hold for a tiling, or a packing of a triangulation, with a set of shapes that are of bounded
Hausdorff distance to circles. At the moment we don't have a proof of the  conjecture even when
the squares are colored black with probability $2/3$.

Behind the second conjecture is a coarse version of conformal invariance. That is, the crossing probability is balanced if the tiles are
of uniformly bounded distance to circles (rotation invariance), and the squares can be of different sizes, (dilation invariance).
\medskip

Let $G$ the $1$-skeleton of bounded degree transient a triangulation of an open disk.
By \cite{HS} it admits a circle packing with similar properties as the tiling in  in conjecture~\ref{c2}.
And if the conformal invariance heuristic holds, we will a.s. see macroscopic crossings for $1/2$-Bernoulli percolation.
\medskip

What about a converse to conjecture 1.1?
\medskip

{\em  Does recurrence implies no percolation
at 1/2?}
\medskip

By same reasoning  we will see unboundely  many macroscopic clusters for $1/2$-Bernoulli percolation, suggesting
that if $G$ is a  $1$-skeleton of bounded degree transient a triangulation of an open disk, then there are
a.s. infinitely many infinite clusters for  $1/2$-Bernoulli site  percolation?
\newpage

{\large 2.1 On the critical probability of planar triangulations. \par}

\medskip

We {\em believe} that $p_c \leq 1/2$  once the triangulation do not have very small (logarithmic) cut sets.

E.g.  if there are $C > 0, \alpha > 0$, so that for every finite set of vertices $S$,

$$
|\partial S| > C |S|^{\alpha}.
$$
Then $p_c \leq 1/2$.

We further believe that $p_c \geq 1/2$ for polynomial growth triangulations of the open disk. Note that if all degrees are at least $6$,
polynomial growth implies that vertices of higher degrees are polynomially sparse, this suggests that their critical probability for  percolation
is $1/2$, as of the triangular lattice. For nonamenable transitive or sofic  triangulations $p_c <1/2$ \cite{BS1}, remove the transitivity assumption.

\medskip

Since we believe that $p_c >0$ for such $G$'s,  by planar duality we conjecture that $p_u < 1$ and {\em uniqueness monotonicity} holds  as well.
Where $p_u$ is the threshold for uniqueness of the infinite cluster.

\section{Conformal invariance and hyperbolicity}

Consider the Poincare disc model of the hyperbolic plane. Pick four points $a,b,c,d$ on the circle at infinity, dividing the circle to four intervals, $A$, $B$, $C$, $D$.

{\em What is the probability that when placing $\lambda$-intensity Poisson process in the disc, with respect to the hyperbolic metric, and coloring
each Voronoi cell black or white independently with equal probability, there is a black crossing between intervals $A$ and $C$ on the boundary?}

Since this process is invariant with respect to hyperbolic isometries, we get that this probability is a function of the crossratio of $a,b,c,d$ and $\lambda$.
There is no scale invariance for the Poisson process on the hyperbolic plane and  increasing $\lambda$  corresponds to the  curvature approaching $0$.

Fix the boundary intervals. It is reasonable to conjecture that the (annealed) crossing probabilities converge as $\lambda$ increases to infinity.
In particular they converge along a subsequence. We get that the subsequential limit is conformal invariant. The limit is Euclidean and should be given by Cardy's formula \cite{Sm}.
The argument above gives the conformal invariance of a subsequential limit.

The point is that conformal invariance of the subsequential limit follows from the hyperbolicity.
Note that in \cite{BS1}, it is shown that  $p_c(\lambda) <1/2$  for any $\lambda$ and suggested that this can be used to show that $p_c \leq 1/2$ in Euclidean Voronoi percolation.
In \cite{BS2}, together  with Oded Schramm  we conjectured that changing the  uniform measure in the disc
(the measure used in sampling the Poisson points) in a uniformly absolutely continuous way, should not effect crossing probability,
as the intensity grows and showed that conformal change of the metric do not effect crossing probabilities.
Here we observe that when placing an infinite measure and unbounded metric,
so that as the intensity grows the local tiling geometry also converges to that of the high intensity Euclidean,
conformality  follows via hyperbolicity.

As in conjecture~\ref{c2} we want: show that the limiting of the crossing probabilities are bounded away from $0$ and $1$ for any non trivial intervals  $A$ and $C$?

\medskip
\noindent
{\bf Acknowledgements:} Thanks to Nicolas Curien, Gady Kozma, Vincent Beffara and Vincent Tassion for useful sanity checks.


\begin{thebibliography}{BKC}

\bibitem{BS0}
I. Benjamini and O. Schramm,
Random walks and harmonic functions on infinite planar graphs using square tilings.
Ann. Probab. 24 (1996) 1219–1238.

\bibitem{BS2}
I. Benjamini and O. Schramm,
Conformal invariance of Voronoi percolation.
Comm. Math. Phys. 197 (1998) 75�-107.


\bibitem{BS1}
I. Benjamini and O. Schramm, Percolation in the hyperbolic plane.
J. Amer. Math. Soc. 14 (2001) 487�-507.


\bibitem{HS}
Z-X. He and O.  Schramm,
Hyperbolic and parabolic packings.
Discrete Comput. Geom. 14 (1995) 123�-149.


\bibitem{Sm}
S. Smirnov, Critical percolation in the plane: conformal invariance, Cardy's formula, scaling limits.
C. R. Acad. Sci. Paris Ser. I Math. 333 (2001)  239-�244.


\bibitem{SW}
S. Smirnov and W.  Werner,
Critical exponents for two-dimensional percolation.
Math. Res. Lett. 8 (2001) 729�-744.


\end{thebibliography}
\end{document}